\documentclass[11pt]{amsart}
\usepackage[colorlinks=true, pdfstartview=FitV, linkcolor=blue, citecolor=blue, urlcolor=blue, breaklinks=true]{hyperref}
\usepackage{amsmath,amsfonts,amssymb,amsthm,amscd,comment,euscript}
\usepackage{xcolor}
\usepackage{float}
\usepackage{gauss} 
\usepackage[latin1]{inputenc}
\usepackage{tikz}
\usepackage[all]{xy}
\usepackage{array}
\usetikzlibrary{calc}

\newcommand{\nc}{\newcommand}

\theoremstyle{remark}
\numberwithin{equation}{section}

\usepackage{amsthm}

\theoremstyle{definition}
\newtheorem{defn}{Definition}[section]
\newtheorem{ex}{Example}[section]
\newtheorem*{cl}{Claim}
\newtheorem*{ack}{Acknowledgements}
\newtheorem{thm}{Theorem}[section]
\newtheorem{rem}{Remark}[section]

\newtheorem{cor}{Corollary}[section]
\newtheorem{lem}{Lemma}[section]
\newtheorem{prop}{Proposition}[section]

\newcommand\Lg{\mathfrak{g}}
\newcommand\U{\mathbf{U}}

\nc{\C}{\mathbb C }
\nc{\Z}{\mathbb Z }
\nc{\N}{\mathbb N }
\newcommand{\wt}{\operatorname{wt}}

\newcommand{\pr}{\operatorname{pr}}
\newcommand{\lev}{\operatorname{lev}}

\begin{document}

\title[Kirillov-Reshetikhin crystals, energy function and the combinatorial R-matrix]{Kirillov-Reshetikhin crystals, energy function and the combinatorial R-matrix}

\author{Deniz Kus}
\address{Deniz Kus:\newline
Mathematisches Institut, Universit\"at zu K\"oln, Germany}
\email{dkus@math.uni-koeln.de}

\thanks{D.K. was partially sponsored by the ``SFB/TR 12-Symmetries and Universality in Mesoscopic Systems".}

\subjclass[2010]{}
\begin{abstract}
We study the polytope model for the affine type $A$ Kirillov-Reshetikhin crystals and prove that the action of the affine Kashiwara operators can be described in a remarkable simple way. Moreover, we investigate the combinatorial $R$-matrix on a tensor product of polytopes and characterize the map explicitly on the highest weight elements. We further give a formula for the local energy function and provide an alternative proof for the perfectness. We determine
 for any dominant highest weight element $\Lambda$ of level $\ell$ the elements $b_{\Lambda}, b^{\Lambda}$ involved in the definition of perfect crystals and  give an explicit description of the ground-state path in the tensor product of polytopes. 
\end{abstract}
\maketitle \thispagestyle{empty}
%%%%%%%%%%%%%%%%%%%%%%%%%%%%%%%%%%%%%%%%%%%%%%%%%%%%%%%%%%%%%%%%%%%%%%%%%%%%%%%%%%%%%%%%%%%%%%%%%%%%%%%%%%%%%%%%%%%%%%%%%%%%%%%%%%%
%         Introduction
%%%%%%%%%%%%%%%%%%%%%%%%%%%%%%%%%%%%%%%%%%%%%%%%%%%%%%%%%%%%%%%%%%%%%%%%%%%%%%%%%%%%%%%%%%%%%%%%%%%%%%%%%%%%%%%%%%%%%%%%%%%%%%%%%%%
\section{Introduction}\label{section1}
Let $\Lg$ be an affine Kac--Moody algebra and $\U^{'}(\Lg)$ be the quantized universal enveloping algebra corresponding to the derived algebra $\Lg^{'}$, called the quantum affine algebra. For finite-dimensional $\U^{'}(\Lg)$-modules $V$ and $V'$, such that the tensor product $V\otimes V^{'}$ is irreducible and $V,V^{'}$ have crystal bases $B$,$B^{'}$ there exists a unique map $R$ from $B\otimes B^{'}$ to $B^{'}\otimes B$ commuting with any Kashiwara operators $\widetilde{e}_l,\widetilde{f}_l$ (see \cite{Ka02}). This map is called the combinatorial $R$-matrix. Moreover, there exists a $\Z$-valued function on $B\otimes B^{'}$ which is defined through a combinatorial rule (see \eqref{locen}) and is called the local energy function. Both functions play an important role in the affine crystal theory. The global energy function is defined (see Definition~\ref{energyd}) on a tensor product $B_1\otimes\cdots\otimes B_N$, where $B_j$ is the crystal basis of a finite-dimensional $\U^{'}(\Lg)$-module $V_j$, through the combinatorial $R$-matrix and the local energy function and is an important grading used in the theory of generalized Kostka polynomials (see \cite{SW99,S02}). The calculation of the combinatorial $R$-matrix or the energy function is done for certain families of crystals in \cite{OS10,SW99,S02}.\par
A certain subclass of finite--dimensional irreducible modules for $\U^{'}(\Lg)$, that gained a lot of attraction during the last decades, are the so-called Kirillov-Reshetikhin modules $W^{i,m}$ where $i$ is a node in the classical Dynkin diagram and $m$ is a positive integer \cite{KR87}. The modules $W^{i,m}$ have distinguished properties among finite-dimensional modules of quantum affine algebras. One of such properties is that Kirillov--Reshetikhin modules were conjectured to admit a crystal bases $B^{i,m}$ (see \cite[Conjecture 2.1]{HKOTT02}) and this was proven for type $A^{(1)}_{n}$ in \cite{KKMMNN92} and for all non-exceptional cases in \cite{OS08}. The combinatorial structure of $B^{i,m}$ was clarified in \cite{FOS09,K12,S02} by exploiting the existence of a map $\sigma$ on $B^{i,m}$ which is the analogue of
the Dynkin diagram automorphism on the level of crystals.\par
Another such property is that KR-crystals were conjectured to be perfect \cite{HKOTT02,HKOTY99}, which is a technical condition, if and only if $m$ is a multiple of a particular constant $c_i$. Perfect crystals are used to give a path realization of crystal bases of integrable highest weight modules for the quantum algebra $\U(\Lg)$ in terms of semi-infinite tensor product of perfect crystals \cite{KKMMNN92}. The highest weight element in this semi-infinite tensor product is called the ground-state path. The prefectness of these crystals was proven for all non-exceptional types in \cite{FOS10}.\par
In this paper we investigate the polytope realization of affine type $A$ Kirillov-Reshetikhin crystals developed in \cite{K12} and consider the natural question of finding explicit formulas for the combinatorial $R$-matrix and the energy function in terms of this realization.\par
Let $A\otimes B$ be a highest weight element in $B^{r_1,s_1}\otimes B^{r_2,s_2}$, then our first result yields the image of $A\otimes B$ under the combinatorial $R$-matrix (see Theorem~\ref{rmatrix}). It is remarkable that the map $R$ behaves almost like the identity map in the sense that it preserves the entries of $A$ and $B$ and changes only the ``shape".\par
Our second result deals with the computation of the energy function. For an arbitrary element $A\otimes B$ contained in $B^{r_1,s_1}\otimes B^{r_2,s_2}$ we associate recursively elements $A\otimes B=A_0\otimes B_0,\ A_1\otimes B_1,\dots,A_k\otimes B_k$, where $k=0$ if $A\otimes B$ is a highest weight element and otherwise depends on the rank of the Lie algebra and $r_2$ (for a precise definition see \eqref{seq}). Our second result shows that the energy of $A\otimes B$ is given up to a sign by the sum over the entries below the $r_2$-th row of $A_k$ (see Theorem~\ref{enonar}).\par
Furthermore, we consider the problem of obtaining the explicit elements $b_{\Lambda}$ and $b^{\Lambda}$ involved in the definition of perfect crystals and the explicit description of the ground-state path. In the language of polytopes these elements can be easily determined and are described in Theorem~\ref{perfectness} where moreover an alternative proof for the perfectness is provided.\par
Another essential part of this paper, which is motivated by a work of Kwon \cite{KW12}, where a combinatorial model for type $A^{(1)}_{n}$ is given in terms of the RSK correspondence, is to simplify the affine crystal structure of $B^{i,m}$ given in \cite{K12} via the promotion operator. Our results give a very explicit and easy way to calculate the classical and affine crystal structure of $B^{i,m}$ even by hand.
The simple formulas for the affine operators, the combinatorial $R$-matrix and the energy function let the author expect that similar calculations can be done for types $B^{(1)}_{n}$ and $C^{(1)}_{n}$ by using the polytopes from \cite{ABS11}. This will be part of forthcoming work.\par
The paper is organized as follows. In Section~\ref{section2} we introduce the main notations and review some general facts about crystals, in particular we recall the realization of affine type $A$ Kirillov-Reshetikhin crystals via polytopes developed in \cite{K12}. In Section~\ref{section3} we recall the notion of Nakajima monomials and prove that the affine crystal structure on the polytope is remarkable simple. In Section~\ref{section4} we determine the classical highest weight elements in the tensor product of polytopes and calculate the image of the combinatorial $R$-matrix on such elements. We also obtain the value of the energy function for arbitrary elements in $B^{r_1,s_1}\otimes B^{r_2,s_2}$. Finally, in Section~\ref{section5} we give an alternative proof for the perfectness of $B^{i,m}$ and determine explicitly the unique elements $b_{\Lambda},b^{\Lambda}$ for any level $\ell$ dominant integral weight $\Lambda$. Further we describe the ground-state path of weight $\Lambda$.
\begin{ack}
The author would like to thank the referee of \cite{K12} for showing strong interest in the descriptions of the combinatorial $R$-matrix and the energy function in terms of this realization.
\end{ack}
%%%%%%%%%%%%%%%%%%%%%%%%%%%%%%%%%%%%%%%%%%%%%%%%%%%%%%%%%%%%%%%%%%%%%%%%%%%%%%%%%%%%%%%%%%%%%%%%%%%%%%%%%%%%%%%%%%%%%%%%%%%%%%%%%%%
%         
%%%%%%%%%%%%%%%%%%%%%%%%%%%%%%%%%%%%%%%%%%%%%%%%%%%%%%%%%%%%%%%%%%%%%%%%%%%%%%%%%%%%%%%%%%%%%%%%%%%%%%%%%%%%%%%%%%%%%%%%%%%%%%%%%%%
\section{Notations and review of crystal theory}\label{section2}
\subsection{Crystal basis and abstract crystals}
Crystal theory provides a combinatorial way to study the representation theory of quantum algebras. In this section we review the theory of crystal bases introduced by Kashiwara in \cite{Kashi91} and fix the main notation. For an indeterminate element $q$ and an affine Lie algebra $\Lg$ with index set $I$ we denote by $\U^{'}_q(\Lg)$ the corresponding quantum algebra without derivation. We denote further by $\U_q(\Lg)$ the corresponding quantum algebra with derivation and by $\U_q(\Lg_0)$ the quantum algebra of the classical subalgebra $\Lg_0$ (with index set $I_0$) of $\Lg$.
A remarkable theorem of Kashiwara implies that every integrable highest weight module $V(\lambda)$ for $\U_q(\Lg)$ (resp. $\U_q(\Lg_0)$) has a crystal basis and hence a corresponding crystal $B(\lambda)$. 
For the finite-dimensional modules of quantum algebras of classical type Kashiwara and Nakashima described $B(\lambda)$ in terms of so-called Kashiwara--Nakashima tableaux \cite{NK94}, the analogue of semi--standard tableaux. For alternative descriptions of $B(\lambda)$ we refer to a series of papers \cite{KKS04,NS97,Ku13,L94}. All these crystal graphs are subject to certain properties, which leads to the definition of abstract crystals.

For us, an abstract crystal is a nonempty set $B$ together with maps
$$\widetilde{e}_l, \widetilde{f}_l: B\longrightarrow B\cup\{0\}, \text{ for } l\in I \text{ (resp. $I_0$)}$$
$$\epsilon_l, \varphi_l: B\longrightarrow \Z, \text{ for } l\in I \text{ (resp. $I_0$) }$$
$$\wt : B\longrightarrow P \text{ (resp. $P_0$) },$$
which satisfy some conditions.
Here $P$ (resp. $P_0$) is the weight lattice associated to $\Lg$ (resp $\Lg_0$). The maps $\widetilde{e}_l$ and $\widetilde{f}_l$ are Kashiwara's crystal operators and $\wt$ is the weight function. 
For quantum algebras of simply-laced Kac--Moody algebras Stembridge gave a set of local axioms characterizing the set of crystals of representations in the class of all crystals \cite{S2003}. Moreover, by using Littelmann's path model Stembridge proved that these axioms hold in all cases, simply-laced or not.

\subsection{Tensor product of crystals and regular crystals}
Suppose that we have two abstract crystals $B_1$, $B_2$, then we can construct a new crystal which is as a set nothing but $B_1\times B_2$. This crystal is denoted by $B_1\otimes B_2$ and the Kashiwara operators are given as follows:
$$\tilde{f}_{l}(b_1\otimes b_2)=\begin{cases} (\tilde{f}_{l}b_1)\otimes b_2, \text{ if $\epsilon_l(b_1)\geq \varphi_l(b_2)$,}\\
b_1\otimes (\tilde{f}_{l}b_2), \text{ if $\epsilon_l(b_1)<\varphi_l(b_2)$.}\end{cases}$$
$$\tilde{e}_{l}(b_1\otimes b_2)=\begin{cases} (\tilde{e}_{l}b_1)\otimes b_2, \text{ if $\epsilon_l(b_1)>\varphi_l(b_2)$,}\\
b_1\otimes (\tilde{e}_{l}b_2), \text{ if $\epsilon_l(b_1)\leq \varphi_l(b_2)$.}\end{cases}$$
Further, one can describe explicitly the maps $\wt$,$\varphi_l$ and $\epsilon_l$ on $B_1\otimes B_2$, namely:
$$\wt(b_1\otimes b_2)=\wt(b_1)+\wt(b_2)$$
$$\varphi_l(b_1\otimes b_2)=\max\{\varphi_l(b_1),\varphi_l(b_1)+\varphi_l(b_2)-\epsilon_l(b_1)\}$$
$$\epsilon_l(b_1\otimes b_2)=\max\{\epsilon_l(b_2),\epsilon_l(b_1)+\epsilon_l(b_2)-\varphi_l(b_2)\}.$$

We say that a crystal $B$ is regular if for each subset $J$ with $|J|=2$ each $J$-component of
$B$ is isomorphic to the crystal of an integrable $\U_q(\Lg_{J})$-module, where $\Lg_J$ is the Kac--Moody
algebra associated to the Cartan matrix $A_J = (a_{i,j})_{i,j\in J}$.

\subsection{Kirillov-Reshetikhin crystals}
The theory of crystal bases can likewise be defined in the setting of $\U^{'}_q(\Lg)$ modules, respecting that crystal bases might not always exist. A certain class of finite-dimensional $\U^{'}_q(\Lg)$ modules, where the existence of crystal bases is proven for the non-exceptional types \cite{OS08}, are the so-called Kirillov--Reshetikhin modules $W^{i,m}$, where $i$ is a node in the classical Dynkin diagram and $m$ is a positive integer. We recall the realization of Kirillov--Reshetikhin crystals for type $A^{(1)}_{n}$ from \cite{K12}. Recall that the classical positive roots are all of the form $$\alpha_{i,j}=\alpha_i+\alpha_{i+1}+\cdots+\alpha_j,\ \mbox{for $1\leq i\leq j\leq n$}.$$

We denote by $\widetilde{B}^{i,m}$ be the set of all pattern

$$
\begin{array}{|c|c|c|c|c|} \hline
a_{1,i} & a_{2,i}& \dots & a_{i-1,i}& a_{i,i}\\[7pt] \hline
a_{1,i+1} & a_{2,i+1} & \dots & \scriptstyle a_{i-1,i+1}& a_{i,i+1} \\[7pt] \hline

\vdots & \vdots & \vdots & \vdots & \vdots \\[2pt] \hline
a_{1,n} & a_{2,n} & \dots &  a_{i-1,n}& a_{i,n} \\[7pt] \hline 
   
\end{array} 
$$
filled with non-negative integers, such that $\sum^n_{s=1}a_{\beta(s)}\leq m$ for all sequences 
$
(\beta(1),\dots, \beta(n))
$
satisfying the following: $\beta(1)=(1,i), \beta(n)=(i,n)$ and if $\beta(s)=(p,q)$ then the next element in the sequence is either of the form 
$\beta(s+1)=(p,q+1)$ or $\beta(s+1)=(p+1,q).$

Further we define the weight function by 
\begin{equation}\label{weight}\wt(A)=m\omega_i-\sum_{ \begin{subarray}{c}1\leq p\leq i\\i\leq q\leq n\end{subarray}}a_{p,q}\alpha_{p,q}.\end{equation}
The maps $\epsilon_l$ and $\varphi_l$ are given as follows
\begin{subequations}
\begin{equation}\label{ttt1}\varphi_l(A)=\begin{cases}
  m-\sum^{i-1}_{j=1} a_{j,i}-\sum^{n}_{j=i} a_{i,j},  & \text{if $l=i$}\\
  \sum^{p^l_+(A)}_{j=1}a_{j,l-1}-\sum^{p^l_+(A)-1}_{j=1}a_{j,l}, & \text{if $l>i$}\\
  \sum^{n}_{j=p^l_{-}(A)}a_{l+1,j}-\sum^{n}_{j=p^l_{-}(A)+1}a_{l,j}, & \text{if $l<i$}
\end{cases}\end{equation}
\begin{equation}\label{ttt2}\epsilon_l(A)=\begin{cases}
  a_{i,i},  & \text{ if $l=i$}\\
  \sum^{i}_{j=q^l_+(A)}a_{j,l}-\sum^{i}_{j=q^l_{+}(A)+1}a_{j,l-1}, & \text{ if $l>i$}\\
  \sum^{q^l_{-}(A)}_{j=i}a_{l,j}-\sum^{q^l_{-}(A)-1}_{j=i}a_{l+1,j}, & \text{ if $l<i$,}
\end{cases}\end{equation}
\end{subequations}
where
\begin{subequations}
\begin{equation}\label{1}p^l_+(A)=\min\bigg\{1\leq p\leq i\mid \sum^p_{j=1}a_{j,l-1}+\sum^{i}_{j=p}a_{j,l}=\max_{1\leq q\leq i}\Big\{\sum^q_{j=1}a_{j,l-1}+\sum^{i}_{j=q}a_{j,l}\Big\}\bigg\}\end{equation}
\begin{equation}\label{2}q^l_+(A)=\max\bigg\{1\leq p\leq i\mid\sum^p_{j=1}a_{j,l-1}+\sum^{i}_{j=p}a_{j,l}=\max_{1\leq q\leq i}\Big\{\sum^q_{j=1}a_{j,l-1}+\sum^{i}_{j=q}a_{j,l}\Big\}\bigg\}\end{equation}
\begin{equation}\label{3}p^l_-(A)=\max\bigg\{i\leq p\leq n\mid\sum^p_{j=i}a_{l,j}+\sum^{n}_{j=p}a_{l+1,j}=\max_{i\leq q\leq n}\Big\{\sum^q_{j=i}a_{l,j}+\sum^{n}_{j=q}a_{l+1,j}\Big\}\bigg\}\end{equation}
\begin{equation}\label{4}q^l_-(A)=\min\bigg\{i\leq p\leq n\mid\sum^p_{j=i}a_{l,j}+\sum^{n}_{j=p}a_{l+1,j}=\max_{i\leq q\leq n}\Big\{\sum^q_{j=i}a_{l,j}+\sum^{n}_{j=q}a_{l+1,j}\Big\}\bigg\}.\end{equation}
\end{subequations}

The Kashiwara operators are defined by 
\begin{subequations}

\begin{equation}\label{kashopf}
\tilde{f}_{l}A=
\begin{cases}
\mbox{replace $a_{i,i}$ by $a_{i,i}+1$}, & \mbox{if $l=i$}\\
\mbox{replace $a_{p^l_+(A),l-1}$ by $a_{p^l_+(A),l-1}-1$ and $a_{p^l_+(A),l}$ by $a_{p^l_+(A),l}+1$},& \mbox{if $l>i$}\\
\mbox{replace $a_{l,p^l_-(A)}$ by $a_{l,p^l_-(A)}+1$ and $a_{l+1,p^l_-(A)}$ by $a_{l+1,p^l_-(A)}-1$}, & \mbox{if $l<i$}\\
\end{cases}
\end{equation}

\begin{equation}\label{kashope}
\tilde{e}_{l}A=
\begin{cases}
\mbox{replace $a_{i,i}$ by $a_{i,i}-1$}, & \mbox{if $l=i$}\\
\mbox{replace $a_{q^l_+(A),l-1}$ by $a_{q^l_+(A),l-1}+1$ and $a_{q^l_+(A),l}$ by $a_{q^l_+(A),l}-1$},& \mbox{if $l>i$}\\
\mbox{replace $a_{l,q^l_-(A)}$ by $a_{l,q^l_-(A)}-1$ and 
$a_{l+1,q^l_-(A)}$ by $a_{l+1,q^l_-(A)}+1$},& \mbox{if $l<i$}.\\

\end{cases}
\end{equation}
\end{subequations}

In order to define the affine operators a map $\pr:\widetilde{B}^{i,m}\longrightarrow \widetilde{B}^{i,m}$ is defined algorithmically (see \cite[Section 5.1]{K12}), which is the analogue of the cyclic Dynkin diagram automorphism $i\mapsto i+1\mod (n+1)$ on the level of crystals. The following theorem gives a realization of the Kirillov--Reshetikhin crystals via polytopes.
\begin{thm}\cite{K12}
The polytope $\widetilde{B}^{i,m}$ with classical crystal structure \eqref{kashopf} and \eqref{kashope} and affine crystal structure 
$$\tilde{f}_0:=\pr^{-1}\circ \tilde{f}_1\circ \pr, \mbox{ and } \tilde{e}_0:=\pr^{-1}\circ  \tilde{e}_1 \circ \pr$$ is isomorphic to the Kirillov--Reshetikhin crystal $B^{i,m}$.
\end{thm}
%%%%%%%%%%%%%%%%%%%%%%%%%%%%%%%%%%%%%%%%%%%%%%%%%%%%%%%%%%%%%%%%%%%%%%%%%%%%%%%%%%%%%%%%%%%%%%%%%%%%%%%%%%%%%%%%%%%%%%%%%%%%%%%%%%%
%         
%%%%%%%%%%%%%%%%%%%%%%%%%%%%%%%%%%%%%%%%%%%%%%%%%%%%%%%%%%%%%%%%%%%%%%%%%%%%%%%%%%%%%%%%%%%%%%%%%%%%%%%%%%%%%%%%%%%%%%%%%%%%%%%%%%%
\section{Simplified affine crystal structure on \texorpdfstring{$B^{i,m}$}{B^{i,m}}}\label{section3}
In this section we define a remarkable simple affine structure on $B^{i,m}$ and keep the explicit classical crystal structure from \eqref{kashopf} and \eqref{kashope}. We prove that $B^{i,m}$ is a regular abstract crystal with respect to the ``new" Kashiwara operators.
\begin{defn}\label{simpleaff}
Let 
$$\varphi_0(A)=a_{1,n} \text{ and } \epsilon_0(A)=m-\sum^n_{j=i}a_{1,j}-\sum^n_{j=2}a_{j,n}.$$
Define $f_0A=0$ (resp. $e_0A=0$) if $\varphi_0(A)=0$ (resp. $\epsilon_0(A)=0$) and otherwise let
$$f_{0}A=\text{ replace }a_{1,n}\text{ by }a_{1,n}-1$$
$$e_{0}A=\text{ replace }a_{1,n}\text{ by }a_{1,n}+1$$
\end{defn}
The maps $f_{0}$ and $e_{0}$ are obviously mutually inverse and we have 

\begin{equation}\label{10}\varphi_0(A)-\epsilon_0(A)=a_{1,n}-m+\sum^n_{j=i}a_{1,j}+\sum^n_{j=2}a_{j,n}=\langle m\omega_i,\alpha^{\vee}_0\rangle+\sum^n_{j=i}a_{1,j}+\sum^n_{j=1}a_{j,n}=\langle\wt(A), \alpha^{\vee}_0\rangle.\end{equation}

\begin{lem}\label{abcry}
The polytope $B^{i,m}$ together with the maps from \eqref{weight}, \eqref{ttt1}, \eqref{ttt2}, \eqref{kashopf},\eqref{kashope} and Definition~\ref{simpleaff} is an abstract crystal.
\proof
The property $\varphi_l(A)-\epsilon_l(A)=\langle \wt(A),\alpha^{\vee}_l\rangle$ for $l=0,\dots,n$ follows from \eqref{10} and \cite[Theorem 3.8.]{K12}. The rest is straightforward. 
\endproof
\end{lem}
\subsection{Nakajima monomials and regularity of \texorpdfstring{$B^{i,m}$}{B^{i,m}}}
In order to prove that $B^{i,m}$ is a regular crystal we shall recall the notion of Nakajima monomials. The set of Nakajima monomials can be understood as a translation of the geometrical realization of crystals provided by Nakajima \cite{Nak01}. Nakajima has shown that there exists a crystal structure on the set of irreducible components of a lagrangian subvariety of the quiver variety. 
%We use this crystal in order to realize the crystal graph $B(\lambda)$, where $\lambda$ is a dominant integral $A_n$ weight, as a set of certain monomials.\par
For $i\in I$ and $n\in\Z$ we consider monomials in the variables $Y_i(n)$, i.e. we obtain the set of Nakajima monomials $\mathcal{M}$ as follows:
$$\mathcal{M}:=\Big\{\prod_{i\in I, n\in\Z} Y_i(n)^{y_i(n)}\mid y_i(n)\in\Z \mbox{ vanish except for finitely many $(i,n)$}\Big\}$$
With the goal to define the crystal structure on $\mathcal{M}$, we take some integers $c=(c_{i,j})_{i\neq j}$ such that $c_{i,j}+c_{j,i}=1$. Let now $M=\prod_{i\in I, n\in\Z} Y_i(n)^{y_i(n)}$ be an arbitrary monomial in $\mathcal{M}$ and $l\in I$, then we set:
$$\wt(M)=\sum_i\big(\sum_n y_i(n)\big)\omega_i$$
$$\varphi_l(M)=\max\Big\{\sum_{k\leq n}y_l(k)\mid n\in\Z\Big\},\ \epsilon_l(M)=\max\Big\{-\sum_{k>n}y_l(k)\mid n\in\Z\Big\}$$
and 
$$n^l_f=\min\Big\{n\mid\varphi_l(M)=\sum_{k\leq n}y_l(k)\Big\},\ n^l_e=\max\Big\{n\mid\epsilon_l(M)=-\sum_{k>n}y_l(k)\Big\}.$$
The Kashiwara operators are defined as follows:
\begin{align*}
&\tilde{f}_{l}M=\begin{cases}
A_l(n^l_f)^{-1}M,& \text{ if $\varphi_l(M)>0$}\\
0,& \text{ if $\varphi_l(M)=0$}\end{cases}&\\&
\tilde{e}_{l}M=\begin{cases}
A_l(n^l_e)M,& \text{ if $\epsilon_l(M)>0$}\\
0,& \text{ if $\epsilon_l(M)=0$,}\end{cases}
\end{align*}
whereby 
$$A_l(n):=Y_l(n)Y_l(n+1)\prod_{i\neq l}Y_i(n+c_{i,l})^{\langle \alpha_i^{\vee},\alpha_l\rangle}.$$

\begin{rem}
\textit{A priori} the crystal structure depends on $c$, hence we will denote this crystal by $\mathcal{M}_c$. But it is easy to see that the isomorphism class of $\mathcal{M}_c$ does not depend on this choice. In the literature $c$ is often chosen as
$c_{ij}=\chi(i\leq j)$ or $c_{ij}=\chi(i\geq j)$ with $\chi$ the indicator function.
\end{rem}

The following result is due to Kashiwara \cite{Kashi03}. 

\begin{prop}\label{iii}
Let $M$ be a monomial in $\mathcal{M}$, such that $\tilde{e}_{l}M=0$ for all $l\in I$. Then the connected component of $\mathcal{M}$ containing $M$ is isomorphic to $B(\wt(M))$.
\end{prop}
The previous proposition will provide a method to prove a simple affine structure on $B^{i,m}$ (cf. Theorem~\ref{thm1} and Corollary~\ref{cor1}). Another method to show this is to verify that the classical crystal structure of \cite{K12} and \cite{KW12} coincide, where the classical crystal structure in \cite{KW12} is obtained from the tensor product of $B(\omega_1)'s$ and is therefore not very explicit. We decided to use the previous proposition, since the classical crystal structures seem quite challenging to compare. 
\begin{thm}\label{thm1}
The polytope $B^{i,m}$ is a regular crystal.
\proof
Let $J=\{k,l\}$ be a subset of $\{0,\dots,n\}$. For $k,l\neq 0$ it is shown in \cite[Theorem 4.5.]{K12} that each $J$-component of $B^{i,m}$ is isomorphic to the crystal of an integrable $\U_q(\Lg_{J})$-module. Hence it remains to prove the statement for $J_1=\{0,1\}$ and $J_n=\{0,n\}$, since the Kashiwara operators $f_0$ and $\widetilde{f}_l$ (resp. $e_0$ and $\widetilde{e}_l$) commute for all $l\neq 1,n$. The proof for $J_1=\{0,1\}$ and $J_n=\{0,n\}$ proceed very similar and therefore we present only the evidence for $J_1=\{0,1\}$. Let $A\in B^{i,m}$ and cancel all arrows in $B^{i,m}$ with color $s\neq 0,1$ and denote the remaining connected graph containing $A$ by $Z_1(A)$. We define a map $\Psi: Z_1(A)\cup \{0\}\longrightarrow \mathcal{M}\cup \{0\}$ which maps $0$ to $0$ and an arbitrary element $B=(b_{p,q})\in Z_1(A)$ to 
$$Y_1(1)^{\sum^n_{j=1}b_{j,n}-m}\prod^{n-i}_{k=0} Y_1(k)^{b_{1,n-k}}\prod^{n-i}_{k=0} Y_2(k)^{b_{2,n-k}}Y_2(k+1)^{-b_{1,n-k}}.$$
By Proposition~\ref{iii} it is enough to prove that the above map $\Psi$ is a strict crystal morphism. We claim 
\begin{equation}\label{gl}n^{2}_f=n-p^1_{-}(B),\ n^{2}_e=n-q^1_{-}(B).\end{equation}
The assumption $p^1_{-}(B)>n-n^{2}_f$ implies 
$$b_{1,n^{2}_f+1}+\cdots+ b_{1,p^1_{-}(B)}\geq b_{2,n^{2}_f}+\cdots+ b_{2,p^1_{-}(B)-1}$$
and thus 
$$\sum^n_{j=n^{2}_f+1}b_{1,j}-\sum^n_{j=n^{2}_f} b_{2,j}\geq \sum^n_{j=p^1_{-}(B)+1}b_{1,j}-\sum^n_{j=p^1_{-}(B)} b_{2,j}.$$
By applying the above inequality by $(-1)$ we arrive at a contradiction to the minimality of $n^{2}_f$.

The assumption $p^1_{-}(B)< n-n^{2}_f$ implies 
$$b_{1,p^1_{-}(B)+1}+\cdots+ b_{1,n^{2}_f}<b_{2,p^1_{-}(B)}+\cdots+ b_{2,n^{2}_f-1}$$
and thus 
$$\sum^n_{j=p^1_{-}(B)+1}b_{1,j}-\sum^n_{j=p^1_{-}(B)}b_{2,j}< \sum^n_{j=n^{2}_f+1}b_{1,j}-\sum^n_{j=n^{2}_f}b_{2,j}.$$
Again by applying the above inequality by $(-1)$ we arrive at a contradiction to the definition of $\varphi_2(\Psi(B))$. Hence $n^{2}_f=n-p^1_{-}(B)$. A similar calculation as above shows $n^{2}_e=n-q^1_{-}(B)$.
Consequently we obtain with \eqref{gl}
\begin{align*}&\varphi_1(B)=\sum^n_{j=p^1_{-}(B)}b_{2,j}-\sum^n_{j=p^1_{-}(B)+1} b_{1,j}=\varphi_2(\Psi(B)),&\\& \epsilon_1(B)=\sum^{q^1_{-}(B)}_{j=i}b_{1,j}-\sum^{q^1_{-}(B)-1}_{j=i} b_{2,j}=\epsilon_2(\Psi(B)).\end{align*}
and if $\widetilde{f}_2(\Psi(B))\neq 0, \widetilde{e}_2(\Psi(B))\neq 0$
\begin{align*}
\widetilde{f}_2(\Psi(B))&=Y_1(1)^{\sum^n_{j=1}b_{j,n}-m}\prod^{n-i}_{k=0} Y_1(k)^{b_{1,n-k}}\prod^{n-i}_{k=0} Y_2(k)^{b_{2,n-k}}Y_2(k+1)^{-b_{1,n-k}}&\\&\hspace{2cm}\times Y_2(n-p^1_{-}(B))^{-1}Y_2(n-p^1_{-}(B)+1)^{-1}Y_1(n-p^1_{-}(B))=\Psi(\widetilde{f}_1B).
\end{align*}
\begin{align*}
\widetilde{e}_2(\Psi(B))&=Y_1(1)^{\sum^n_{j=1}b_{j,n}-m}\prod^{n-i}_{k=0} Y_1(k)^{b_{1,n-k}}\prod^{n-i}_{k=0} Y_2(k)^{b_{2,n-k}}Y_2(k+1)^{-b_{1,n-k}}&\\&\hspace{2cm}\times Y_2(n-q^1_{-}(B))Y_2(n-q^1_{-}(B)+1)Y_1(n-q^1_{-}(B))^{-1}=\Psi(\widetilde{e}_1B).
\end{align*}
For $l=1$ we get
\begin{align*}\varphi_1(\Psi(B))=\max\big\{b_{1,n},b_{1,n}+\underbrace{\sum^n_{j=1}b_{j,n}-m}_{\leq0},b_{1,n}+\underbrace{b_{1,n-1}+\sum^n_{j=1}b_{j,n}-m}_{\leq0},\dots\big\}=b_{1,n},\quad n^1_{f}=0\end{align*}
and 
\begin{align*}\epsilon_1(\Psi(B))=\max\big\{-b_{1,i},-b_{1,i}-b_{1,i+1},\dots,-\sum^{n-2}_{j=i}b_{1,j},m-\sum^{n-1}_{j=i}b_{1,j}-\sum^n_{j=1}b_{j,n}\big\}=\epsilon_0(B),\quad n^1_{e}=0.\end{align*}
Therefore, when $\widetilde{f}_1\Psi(B)\neq 0, \widetilde{e}_1\Psi(B)\neq 0$ we obtain
\begin{align*}\widetilde{f}_1\Psi(B)&=Y_1(1)^{\sum^n_{j=1}b_{j,n}-m}\prod^{n-i}_{k=0} Y_1(k)^{b_{1,n-k}}\prod^{n-i}_{k=0} Y_2(k)^{b_{2,n-k}}Y_2(k+1)^{-b_{1,n-k}}&\\&\hspace{2cm}\times Y_1(0)^{-1}Y_1(1)^{-1}Y_2(1)=\Psi(f_0B)\end{align*}
\begin{align*}\widetilde{e}_1\Psi(B)&=Y_1(1)^{\sum^n_{j=1}b_{j,n}-m}\prod^{n-i}_{k=0} Y_1(k)^{b_{1,n-k}}\prod^{n-i}_{k=0} Y_2(k)^{b_{2,n-k}}Y_2(k+1)^{-b_{1,n-k}}&\\&\hspace{2cm}\times Y_1(0)Y_1(1)Y_2(1)^{-1}=\Psi(e_0B),\end{align*}
which finishes the proof.
\endproof
\end{thm}
As a corollary we obtain a remarkable simple affine crystal structure and an explicit classical crystal structure on the Kirillov-Reshetikhin crystal $B^{i,m}$.
\begin{cor}\label{cor1}
We have 
$$\widetilde{f}_0=\pr^{-1}\circ \widetilde{f}_0\circ \pr=f_0,\ \widetilde{e}_0=\pr^{-1}\circ \widetilde{e}_0\circ \pr=e_0.$$
\proof
The proof follows from Theorem~\ref{thm1}, \cite[Theorem 4.5.]{K12} and \cite[Lemma 2.6]{ST12}.
\endproof
\end{cor}

%%%%%%%%%%%%%%%%%%%%%%%%%%%%%%%%%%%%%%%%%%%%%%%%%%%%%%%%%%%%%%%%%%%%%%%%%%%%%%%%%%%%%%%%%%%%%%%%%%%%%%%%%%%%%%%%%%%%%%%%%%%%%%%%%%%
%         
%%%%%%%%%%%%%%%%%%%%%%%%%%%%%%%%%%%%%%%%%%%%%%%%%%%%%%%%%%%%%%%%%%%%%%%%%%%%%%%%%%%%%%%%%%%%%%%%%%%%%%%%%%%%%%%%%%%%%%%%%%%%%%%%%%%
\section{Combinatorial R-matrix and the energy function}\label{section4}
\subsection{Combinatorial R-matrix}
Let $B_1$ and $B_2$ be two affine crystals with generators $v_1$ and $v_2$ such that the tensor product $B_1\otimes B_2$ is connected and $v_1\otimes v_2$ lies in a one-dimensional weight space. The combinatorial R-matrix (see \cite[Section 4]{KKMMNN1992}) is the unique affine crystal isomorphism 
$$\sigma: B_{1}\otimes B_{2}\stackrel{\sim}{\longrightarrow} B_{2}\otimes B_{1}.$$
We consider two Kirillov-Reshetikhin crystals $B^{r_1,s_1}$,$B^{r_2,s_2}$ from Section~\ref{section2} with the simplified crystal structure proven in Section~\ref{section3}.
By weight consideration we must have $\sigma(v_1\otimes v_2)=v_2\otimes v_1,$ where the generator $v_j\in B^{r_j,s_j}$ is the unique element with zero entries.
In the following we determine the combinatorial $R$-matrix on the classical highest weight vectors of $B^{r_1,s_1}\otimes B^{r_2,s_2}$. The highest weight elements are described in the following lemma. Let $\mathbf s=\min\{s_1,s_2\}$, $\mathbf r=\min\{r_1,r_2\}$, $\mathbf{\widetilde{r}}=\max\{r_1,r_2\}$ and $\mathbf k=\min \big\{\mathbf r-1,n-\mathbf{\widetilde{r}}\big\}$.
\begin{lem}\label{hwe}
The set of highest weight elements is in bijection to 
$$\Big\{(a_{0},a_{1},\dots,a_{k})\in \Z_{\geq 0}^{\mathbf{k}+1}\mid 0\leq a_{\mathbf{k}}\leq a_{\mathbf{k}-1}\leq \cdots\leq a_{0}\leq \mathbf s\Big\}.$$
To be more precise, if $A\otimes B\in B^{r_1,s_1}\otimes B^{r_2,s_2}$ is a highest weight element then we have
$$B=0 \text{ and } a_{r,s}=0 \text{ for all } (r,s)\notin \big\{(\mathbf r,\mathbf{\widetilde{r}}),(\mathbf r-1,\mathbf{\widetilde{r}}+1),\dots,(\mathbf r-\mathbf{k},\mathbf{\widetilde{r}}+\mathbf{k})\big\}$$
\text{and} 
$$0\leq a_{\mathbf r-\mathbf{k},\mathbf{\widetilde{r}}+\mathbf{k}}\leq\dots\leq a_{\mathbf r,\mathbf{\widetilde{r}}}\leq \mathbf s.$$ 
\end{lem}
%\begin{rem}\label{hwe2}
%If $r_1\geq r_2$ and $A\otimes B\in B^{r_1,s_1}\otimes B^{r_2,s_2}$ is a highest weight element, then we have
%$$B=0 \text{ and } a_{r,s}=0 \text{ for all } (r,s)\notin \big\{(r_2,r_1),(r_2-1,r_1+1),\dots,(r_2-k,r_1+k)\big\}$$
%\text{and} 
%$$0\leq a_{r_2-k,r_1+k}\leq \cdots\leq a_{r_2,r_1}\leq \min\{s_1,s_2\},$$ 
%$$\Big\{(a_{r_2,r_1},a_{r_2-1,r_1+1},\dots,a_{r_2-k,r_1+k})\in \Z_{\geq0}^{k+1}\mid 0\leq a_{r_2-k,r_1+k}\leq \cdots\leq a_{r_2,r_1}\leq \min\{s_1,s_2\}\Big\},$$
%where $k=\min\{r_2-1,n-r_1\}.$
%\end{rem}
\proof
We prove the lemma for $r_1\leq r_2$, since the case $r_1\leq r_2$ proceeds similarly. 
The tensor product property from Section~\ref{section2} and $\widetilde{e}_l(A\otimes B)=0$ implies
$$\widetilde{e}_lB=0 \text{ for all } l=1,\dots,n \text{ and } \widetilde{e}_{l}A=0 \text{ for all } l\neq r_2 \text{ and } \epsilon_{r_2}(A)\leq s_2.$$
Hence $B=0$. We claim the following:\\
\begin{cl}
For $A\in B^{r_1,s_1}$ with $\epsilon_j(A)=0$ for all $j\neq r_2$ we have 
$$a_{r,s}=0 \text{ for all } (r,s)\notin \big\{(r_1,r_2),(r_1-1,r_2+1),\dots,(r_1-\mathbf{k},r_2+\mathbf{k})\big\}.$$
\end{cl}
We proof the claim by upward induction on $n$, where the initial step ($n=1$) is obviously true. The property $\epsilon_j(A)=0$ for all $j\neq r_2$ requires 
$$a_{1,r_1}=\cdots=a_{r_1,r_1}=a_{r_1,r_1+1}=\cdots=\widehat{a_{r_1,r_2}}=\cdots=a_{r_1,n}=0.$$ 
Moreover, it is easy to see that $a_{r,s}=0$ for all $r_1\leq s\leq r_2$ and $1\leq r\leq r_1-1$. Hence the $r_1=1$ case is done, so let $r_1>1$. We denote by $\widetilde{A}$ the element obtained from $A$ by removing the last column. We consider $\widetilde{A}$ as an element in $B^{r_1-1,s_1}$ for $A^{(1)}_{n-1}$. We shall prove
\begin{equation}\label{claim1}\epsilon_j(\widetilde{A})=0 \text{ for all } j\neq r_2+1, \end{equation}
where $\epsilon_{r_2}(\widetilde{A})=0$ follows from the above observation.
Assume $\epsilon_j(\widetilde{A})>0$ for some $j\neq r_2,r_2+1$ and $j\geq r_1$. Since $j\neq r_2+1$ we have $a_{r_1,j-1}=0$ and thus 
$$0<\epsilon_j(\widetilde{A})\leq \epsilon_j(\widetilde{A})+a_{r_1,j}=\epsilon_j(A),$$
which is a contradiction to $\epsilon_j(A)=0$. Now assume $\epsilon_j(\widetilde{A})>0$ for some $j$ with $j<r_1$. Then 
$$0<\epsilon_j(\widetilde{A})=\epsilon_j(\widetilde{A})+a_{j,r_1}-a_{j+1,r_1}=\epsilon_j(A)=0,$$
which is once more a contradiction and therefore we get \eqref{claim1}. By induction we obtain that $a_{r,s}=0 \text{ for all } (r,s)\notin \big\{(r_1-1,r_2+1),\dots,(r_1-\mathbf{k},r_2+\mathbf{k})\big\},$
which finishes the proof of the claim.\par
Since $A\in B^{r_1,s_1}$ and $\epsilon_{r_2}(A)=a_{r_1,r_2}\leq s_2$ we get $a_{r_1,r_2}\leq \mathbf s.$ The property $a_{r_1-s,r_2+s}<a_{r_1-s-1,r_2+s+1}$ would imply $\epsilon_{r_2+s+1}(A)> 0$ and thus $0\leq a_{r_1-\mathbf{k},r_2+\mathbf{k}}\leq \cdots\leq a_{r_1,r_2}\leq \mathbf s$.
Now it remains to prove that all elements with the property from Lemma~\ref{hwe} are highest weight elements, which is straightforward.\endproof
The next theorem describes the image of $\sigma$ restricted to the set of highest weight vectors in $B^{r_1,s_1}\otimes B^{r_2,s_2}$. The isomorphism $\sigma$ on arbitrary elements is therefore reduced to a classical problem.
\begin{thm}\label{rmatrix}
Let $A\otimes B$ be a classical highest weight element in $B^{r_1,s_1}\otimes B^{r_2,s_2}$. Then $\sigma(A\otimes B)=\widetilde{A}\otimes \widetilde{B}$ is the unique highest weight element in $B^{r_2,s_2}\otimes B^{r_1,s_1}$ with 
$$\widetilde{B}=0,\text{ and }\widetilde{a}_{\mathbf{r}-j,\mathbf{\widetilde{r}}+j}=a_{\mathbf{r}-j,\mathbf{\widetilde{r}}+j} \text{ for all } j=0,1,\dots,\mathbf{k}.$$
\proof
Again we provide the proof only for $r_1\leq r_2$. The combinatorial R-matrix is a crystal isomorphism and thus commutes with all operators $\widetilde{e}_l,\widetilde{f}_l$. This implies that $\sigma(A\otimes B)$ is again a classical highest weight element. By Lemma~\ref{hwe} we have
$$\widetilde{a}_{r,s}=0, \text{ for all } (r,s)\notin\big\{(r_1,r_2),(r_1-1,r_2+1),\dots,(r_1-\mathbf{k},r_2+\mathbf{k})\big\}$$
and $0\leq \widetilde{a}_{r_1-\mathbf{k},r_2+\mathbf{k}}\leq \cdots\leq \widetilde{a}_{r_1,r_2}\leq \mathbf s.$ We prove by downward induction on $j$ that $\widetilde{a}_{r_1-j,r_2+j}=a_{r_1-j,r_2+j}$ for all $j=\mathbf{k},\mathbf{k}-1,\dots,0$. 
We start with $j=\mathbf{k}$ and suppose $\mathbf{k}=n-r_2\leq r_1-1$. Since $r_1+r_2-n-1\neq r_1,r_2$ we have
\begin{align*}\varphi_{r_1+r_2-n-1}(A\otimes 0)&=\max\big\{a_{r_1+r_2-n,n}, a_{r_1+r_2-n,n}+\varphi_{r_1+r_2-n-1}(0)-\epsilon_{r_1+r_2-n-1}(A)\big\}&\\&=\max\big\{a_{r_1+r_2-n,n}, a_{r_1+r_2-n,n}+0-0\big\}=a_{r_1+r_2-n,n}&\\&=\varphi_{r_1+r_2-n-1}(\widetilde{A}\otimes 0)&\\&=\max\big\{\widetilde{a}_{r_1+r_2-n,n}, \widetilde{a}_{r_1+r_2-n,n}+\varphi_{r_1+r_2-n-1}(0)-\epsilon_{r_1+r_2-n-1}(\widetilde{A})\big\}=\widetilde{a}_{r_1+r_2-n,n}.\end{align*} For $k=r_1-1<n-r_2$ we get similarly as above with $r_1+r_2\neq r_1,r_2$
\begin{align*}\varphi_{r_1+r_2}(A\otimes 0)&=\max\big\{a_{1,r_1+r_2-1}, a_{1,r_1+r_2-1}+\varphi_{r_1+r_2}(0)-\epsilon_{r_1+r_2}(A)\big\}&\\&=\max\big\{a_{1,r_1+r_2-1}, a_{1,r_1+r_2-1}+0-0\big\}=a_{1,r_1+r_2-1}&\\&=\varphi_{r_1+r_2}(\widetilde{A}\otimes 0)&\\&=\max\big\{\widetilde{a}_{1,r_1+r_2-1}, \widetilde{a}_{1,r_1+r_2-1}+\varphi_{r_1+r_2}(0)-\epsilon_{r_1+r_2}(\widetilde{A})\big\}=\widetilde{a}_{1,r_1+r_2-1},\end{align*}
which finishes the initial step. 
Now let $j<\mathbf{k}$. Since $r_2+j+1\neq r_1,r_2$ (otherwise $j+1=0$ or $j+1=r_1-r_2\leq 0$) we obtain by using the induction hypothesis 
\begin{align*}\varphi_{r_2+j+1}(A\otimes 0)&=\max\big\{a_{r_1-j,r_2+j}-a_{r_1-j-1,r_2+j+1},\varphi_{r_2+j+1}(A)+\varphi_{r_2+j+1}(0)-\epsilon_{r_2+j+1}(A)\big\}&\\&=\max\big\{a_{r_1-j,r_2+j}-a_{r_1-j-1,r_2+j+1},a_{r_1-j,r_2+j}-a_{r_1-j-1,r_2+j+1}+0-0\big\}&\\&=a_{r_1-j,r_2+j}-a_{r_1-j-1,r_2+j+1}&\\&=\varphi_{r_2+j+1}(\widetilde{A}\otimes 0)&\\&=\max\big\{\widetilde{a}_{r_1-j,r_2+j}-\widetilde{a}_{r_1-j-1,r_2+j+1},\varphi_{r_2+j+1}(\widetilde{A})+\varphi_{r_2+j+1}(0)-\epsilon_{r_2+j+1}(\widetilde{A})\big\}&\\&=\widetilde{a}_{r_1-j,r_2+j}-\widetilde{a}_{r_1-j-1,r_2+j+1}=\widetilde{a}_{r_1-j,r_2+j}-a_{r_1-j-1,r_2+j+1},\end{align*}
which finishes the proof of the theorem.
\endproof
\end{thm}
We have seen that the map $\sigma$ behaves on the classical highest weight elements almost like the identity map. To find a formula for $\sigma$ on arbitrary elements seems to be quite challenging, which is illustrated in the following example.
\begin{ex}
We consider the Lie algebra $A^{(1)}_{1}$ and the tensor product of Kirillov-Reshetikhin crystals $B^{1,s_1}\otimes B^{1,s_2}$ with $s_1\leq s_2$. For $A\otimes B\in B^{1,s_1}\otimes B^{1,s_2}$ we denote by $\sigma(A)\otimes \sigma(B)$ the image of $A\otimes B$ in 
$B^{1,s_2}\otimes B^{1,s_1}$. The study of the maps $\varphi_l$ and a simple weight consideration determines $\sigma(A),\sigma(B)$. We shall consider three cases:
\begin{align*}
&\bullet \text{ Let }A+B\leq s_1,\text{ then }
\sigma(A)=A,\ \sigma(B)=B.&\\&
\bullet \text{ Let }s_1<A+B\leq s_2,\text{ then }
\sigma(A)=2A-s_1+B,\ \sigma(B)=s_1-A.&\\&
\bullet \text{ Let }s_2<A+B,\text{ then }
\sigma(A)=A+s_2-s_1,\ \sigma(B)=s_1-s_2+B.\end{align*}

Below are the crystal graphs of $B^{1,1}\otimes B^{1,3}$ and $B^{1,3}\otimes B^{1,1}$ respectively 

$$
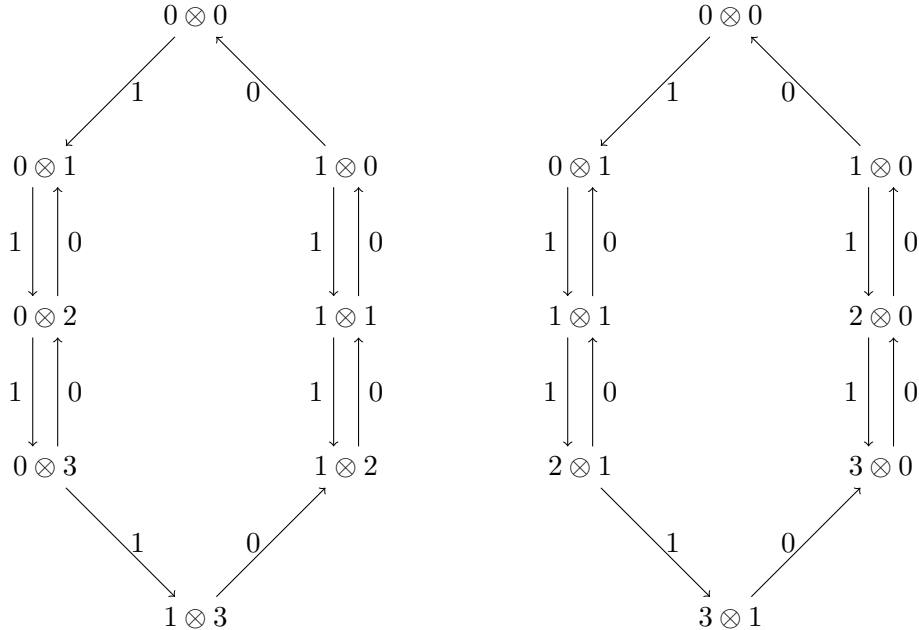
\begin{tikzpicture}

  \node (n1) at (-7,5)  {$0\otimes 0$};
  \node (n2) at (-9,3)  {$0\otimes 1$};
  \node (n3) at (-9,1)  {$0\otimes 2$};
  \node (n4) at (-9,-1)  {$0\otimes 3$};
  \node (n5) at (-5,-1)  {$1\otimes 2$};
  \node (n6) at (-5,1)  {$1\otimes 1$};
  \node (n7) at (-5,3)  {$1\otimes 0$};
  \node (n8) at (-7,-3)  {$1\otimes 3$};

  \draw[->] (n1)edge node[right,midway] {1}(n2);
  \draw[transform canvas={xshift=-1ex},->] (n2)edge node[left,midway] {1}(n3);
  \draw[transform canvas={xshift=-1ex},->] (n3)edge node[left,midway] {1}(n4);
  \draw[->] (n4)edge node[right,midway] {1}(n8);
  \draw[transform canvas={xshift=-1ex},->] (n7)edge node[left,midway] {1}(n6);
  \draw[transform canvas={xshift=-1ex},->] (n6)edge node[left,midway] {1}(n5);
  \draw[->] (n8)edge node[left,midway] {0} (n5);
  \draw[transform canvas={xshift=1ex},->] (n5)edge node[right,midway] {0} (n6);
  \draw[transform canvas={xshift=1ex},->] (n6)edge node[right,midway] {0} (n7);
  \draw[->] (n7)edge node[left,midway] {0} (n1);
  \draw[transform canvas={xshift=1ex},->] (n4)edge node[right,midway] {0} (n3);
  \draw[transform canvas={xshift=1ex},->] (n3)edge node[right,midway] {0} (n2);
  \end{tikzpicture}
\hspace{2cm}
\begin{tikzpicture}

  \node (n1) at (-7,5)  {$0\otimes 0$};
  \node (n2) at (-9,3)  {$0\otimes 1$};
  \node (n3) at (-9,1)  {$1\otimes 1$};
  \node (n4) at (-9,-1)  {$2\otimes 1$};
  \node (n5) at (-5,-1)  {$3\otimes 0$};
  \node (n6) at (-5,1)  {$2\otimes 0$};
  \node (n7) at (-5,3)  {$1\otimes 0$};
  \node (n8) at (-7,-3)  {$3\otimes 1$};

  \draw[->] (n1)edge node[right,midway] {1}(n2);
  \draw[transform canvas={xshift=-1ex},->] (n2)edge node[left,midway] {1}(n3);
  \draw[transform canvas={xshift=-1ex},->] (n3)edge node[left,midway] {1}(n4);
  \draw[->] (n4)edge node[right,midway] {1}(n8);
  \draw[transform canvas={xshift=-1ex},->] (n7)edge node[left,midway] {1}(n6);
  \draw[transform canvas={xshift=-1ex},->] (n6)edge node[left,midway] {1}(n5);
  \draw[->] (n8)edge node[left,midway] {0} (n5);
  \draw[transform canvas={xshift=1ex},->] (n5)edge node[right,midway] {0} (n6);
  \draw[transform canvas={xshift=1ex},->] (n6)edge node[right,midway] {0} (n7);
  \draw[->] (n7)edge node[left,midway] {0} (n1);
  \draw[transform canvas={xshift=1ex},->] (n4)edge node[right,midway] {0} (n3);
  \draw[transform canvas={xshift=1ex},->] (n3)edge node[right,midway] {0} (n2);
\end{tikzpicture}
$$
The study of the map $\sigma$ on arbitrary elements for higher rank Lie algebras shows that we need considerably more case differentiations. 
\end{ex}
\begin{ex}
Consider the Lie algebra $A^{(1)}_{7}$ and the affine crystal $B^{4,s_1}\otimes B^{5,s_2}$. Then 
$$
\begin{array}{|c|c|c|c|c|} \hline
0 & 0 & 0& 0\\ \hline
0 & 0 & 0& a \\ \hline
0 & 0 & b& 0 \\ \hline
0 & c & 0& 0 \\ \hline
\end{array} 
\otimes
\begin{array}{|c|c|c|c|c|} \hline
0 & 0 & 0& 0& 0\\ \hline
0 & 0 & 0& 0& 0 \\ \hline
0 & 0 & 0& 0& 0 \\ \hline
\end{array} 
\stackrel{\sigma}{\longrightarrow}
\begin{array}{|c|c|c|c|c|} \hline
0 & 0 & 0& a& 0\\ \hline
0 & 0 & b& 0& 0 \\ \hline
0 & c & 0& 0& 0 \\ \hline
\end{array}
\otimes
\begin{array}{|c|c|c|c|c|} \hline
0 & 0 & 0& 0\\ \hline
0 & 0 & 0& 0 \\ \hline
0 & 0 & 0& 0 \\ \hline
0 & 0 & 0& 0 \\ \hline
\end{array} 
 $$

\end{ex}
\subsection{Energy function}
We recall the definition of the energy function from \cite{KKMMNN1992,ST12}. There exists a function called the local energy function $H=H_{B_1,B_2}: B_1\otimes B_2\longrightarrow \Z$ unique up to global additive constant, such that
\begin{equation}\label{locen}H\big(e_l(b_1\otimes b_2)\big)=H(b_1\otimes b_2)+\begin{cases} -1& \text{ if $l=0$ and LL}\\
1& \text{ if $l=0$ and RR}\\
0& \text{ otherwise}\end{cases},\end{equation}
where LL (resp. RR) means that $\widetilde{e}_0$ acts on both $b_1\otimes b_2$ and $\sigma(b_1\otimes b_2)$ on the first (resp. second) tensor factor.
We consider the tensor product of KR-crystals and normalize $H$ by requiring $H(0\otimes0)=0$.
\begin{defn}\label{energyd}
For $B=B^{r_1,s_2}\otimes\cdots\otimes B^{r_N,s_N}$ and $1\leq i<j\leq N$ set
$$H_{j,i}=H_i\sigma_{i+1}\cdots\sigma_{j-1},$$
where $\sigma_i$ and $H_i$ act on the $i$-th and $(i+1)$-th tensor factor. The energy function is defined as
$$D_B=\sum_{1\leq i<j\leq N} H_{i,j}.$$
\end{defn}
\begin{rem}
In the definition of the energy function for tensor products of KR-crystals of arbitrary type the so-called $D$-function $D:B^{r_i,s_i}\longrightarrow \Z$ is involved, which is constant on all classical components. Since the KR-crystals for type $A^{(1)}_{n}$ are classicaly irreducible $D$ is the constant function 0.
\end{rem}
The aim of the rest of this section is to give an explicit formula for the local energy function. Together with this formula and Theorem~\ref{rmatrix} we obtain a formula for the energy function $D_B$. 
\begin{prop}\label{energy1}
Let $A\otimes B$ be a highest weight element in $B^{r_1,s_1}\otimes B^{r_2,s_2}$. 
Then 
$$H(A\otimes B)=-\sum_{\begin{subarray}{c}1\leq p\leq r_1\\r_1\leq q\leq n\end{subarray}}a_{p,q}.$$
In particular, the local energy function is the negative sum over the entries of $A$.
\proof
We assume $r_1\leq r_2$. We know by Lemma~\ref{hwe} that the only possible non-zero entries of $A$ are given by
$$a_{r_1,r_2},a_{r_1-1,r_2+1},\dots,a_{r_1-\mathbf{k},r_2+\mathbf{k}}.$$
We prove the statement by induction on $\sum^{\mathbf{k}}_{j=0}a_{r_1-j,r_2+j}$. When $\sum^{\mathbf{k}}_{j=0}a_{r_1-j,r_2+j}=0$ we have by normalization $H(0\otimes 0)=0$ and the statement is true. So suppose that the sum is greater than zero and let $j$ be maximal with $\mathbf a=a_{r_1-j,r_2+j}\neq 0$. In the following we exploit the simplified crystal structure from Section~\ref{section3}. It is easy to see that there exists a tuple $(i_1,\dots,i_p)$ with $i_p\notin\{0,r_1,r_2\}$ and $\varphi_0(\widetilde{f}_{i_p}\cdots \widetilde{f}_{i_1}A)\geq \mathbf a$. For instance we can choose the tuple  
$$\big((r_1-j-1)^{\mathbf a},(r_1-j-2)^{\mathbf a},\dots,1^{\mathbf a},(r_2+j+1)^{\mathbf a},(r_2+j+2)^{\mathbf a},\dots,n^{\mathbf a}\big).$$
Let $\widetilde{A}=\widetilde{f}_{0}^{\mathbf a}\widetilde{f}_{i_p}\cdots \widetilde{f}_{i_1}A$, which is the element obtained from $A$ by replacing $\mathbf a$ by $0$. Especially $\widetilde{A}\otimes 0$ is a highest weight element and by induction we get
$$H(\widetilde{A}\otimes0)=-\sum_{\begin{subarray}{c}1\leq p\leq r_1\\r_1\leq q\leq n\end{subarray}}a_{p,q}+\mathbf a.$$
Consequently, we obtain (recall $i_p\notin\{0,r_1,r_2\}$)
$$H(A\otimes0)=H\big(\widetilde{e}^{\mathbf a}_0\widetilde{f}^{\mathbf a}_0\widetilde{f}_{i_p}\cdots \widetilde{f}_{i_1}(A\otimes0)\big)=H\big(\widetilde{e}^{\mathbf a}_0(\widetilde{A}\otimes 0)\big)=H\big(\widetilde{e}^{\mathbf{a}-1}_0(\widetilde{A}\otimes0)\big)+\begin{cases} -1& \text{ if  LL}\\
1& \text{ if  RR}\\
0& \text{ otherwise.}\end{cases}$$
In order to prove the proposition we shall show that $\widetilde{e}_0$ acts on $\widetilde{e}^{\mathbf{a}-t}_0(\widetilde{A}\otimes0)$ and $\sigma\big(\widetilde{e}^{\mathbf{a}-t}_0(\widetilde{A}\otimes0)\big)$ on the first tensor factor for all $1\leq t\leq \mathbf a$. Assume we have that the second tensor factor in $\sigma\big(\widetilde{e}^{\mathbf{a}-t}_0(\widetilde{A}\otimes0)\big)$ is zero, i.e.
\begin{equation}\label{claim2}\sigma\big(\widetilde{e}^{\mathbf{a}-t}_0(\widetilde{A}\otimes0)\big)=\sigma\big(\widetilde{e}^{\mathbf{a}-t}_0\widetilde{A}\big)\otimes 0.\end{equation}
Then we get
$$\epsilon_0(\widetilde{e}^{\mathbf{a}-t}_0\widetilde{A})=s_1-\mathbf{a}+t>\varphi_0(0)=0$$
and
$$\epsilon_0\big(\sigma\big(\widetilde{e}^{\mathbf{a}-t}_0\widetilde{A}\big)\otimes 0\big)=\epsilon_0\big(\sigma\big(\widetilde{e}^{\mathbf{a}-t}_0\widetilde{A}\big)\big)+s_1=\epsilon_0\big(\widetilde{e}^{\mathbf{a}-t}_0(\widetilde{A}\otimes 0)\big)=s_1-\mathbf a+t+s_2.$$

$$\Longrightarrow \epsilon_0\big(\sigma\big(\widetilde{e}^{\mathbf{a}-t}_0\widetilde{A}\big)\big)=s_2-\mathbf a+t>\varphi_0(0)=0,$$
which gives the statement. Hence it remains to prove \eqref{claim2}.
We set
$\widetilde{\widetilde{A}}\otimes 0=\widetilde{e}_{j_r}\cdots\widetilde{e}_{j_1}\widetilde{e}^{\mathbf a-t}_0(\widetilde{A}\otimes 0),$ where $(j_1,\dots,j_r)$ is the tuple
$$\big(n^{\mathbf a-t},(n-1)^{\mathbf a-t},\dots,(r_2+j+1)^{\mathbf a-t},1^{\mathbf a-t},\dots,(r_1-j-2)^{\mathbf a-t},(r_1-j-1)^{\mathbf a-t}\big).$$
In particular $\widetilde{\widetilde{A}}\otimes 0$ is a highest weight element and $\widetilde{\widetilde{A}}$ is obtained from $A$ by replacing $\mathbf a$ by $\mathbf a-t$. By Theorem~\ref{rmatrix} we obtain that $\sigma(\widetilde{\widetilde{A}}\otimes 0)$ is of the form $\sigma(\widetilde{\widetilde{A}})\otimes 0$.
Furthermore,
$$\sigma\big(\widetilde{e}^{\mathbf{a}-t}_0(\widetilde{A}\otimes0)\big)=\widetilde{f}_{j_1}\cdots\widetilde{f}_{j_r}\sigma(\widetilde{\widetilde{A}}\otimes 0).$$
Thus $\sigma\big(\widetilde{e}^{\mathbf{a}-t}_0(\widetilde{A}\otimes0)\big)$ is of the form $\sigma\big(\widetilde{e}^{\mathbf{a}-t}_0\widetilde{A}\big)\otimes 0,$ since $j_p\notin\{r_1,r_2\}$ for all $p\in\{1,\dots,r\}.$
The proof for $r_1\geq r_2$ proceeds similarly.
\endproof
\end{prop}
In the remaining part of this subsection we determine the local energy function on arbitrary elements $A\otimes B\in B^{r_1,s_1}\otimes B^{r_2,s_2}$. We define a map $+:\Z\longrightarrow\Z_{\geq 0},\ x\mapsto x_+$, where $x_+=x$, if $x>0$ and $x_+=0$ otherwise.\par
Further we define elements $A_{r}^{s}=\big((a^r_s)_{p,q}\big)$, $B_{r}^{s}=\big((b^r_s)_{p,q}\big)$ for $0\leq s\leq n-r_2$ and $0\leq r\leq r_2+s$ as follows:
\begin{align}\label{seq}
\text{For $s=0$:}\hspace{1cm}&A_{0}^{0}=A,\hspace{5,02cm} B_{0}^{0}=B&\\& \nonumber  A_{r}^{0}=\widetilde{e}^{\big(\epsilon_r(A^{0}_{r-1})-\varphi_r(B^{0}_{r-1})\big)_{+}}_rA^{0}_{r-1},\hspace{1cm} B^{0}_{r}=\widetilde{e}_r^{\epsilon_r(B^{0}_{r-1})}B^{0}_{r-1}.\end{align}
$$\hspace{0,7cm}\text{For $s>0$:} \hspace{1cm}A_{r}^{s}=\begin{cases}A_{r_2+s-1}^{s-1} &\text{if $r=0$}\\
\widetilde{e}^{\big(\epsilon_r(A^{s}_{r-1})-\varphi_r(B^{s}_{r-1})\big)_{+}}_rA_{r-1}^{s} &\text{if $1\leq r\leq r_2-1$}\\
\widetilde{e}^{\big(\epsilon_{2r_2+s-r}(A^{s}_{r-1})-\varphi_{2r_2+s-r}(B^{s}_{r-1})\big)_{+}}_{2r_2+s-r}A_{r-1}^{s} &\text{if $r_2\leq r\leq r_2+s$}\end{cases}
$$
$$
B_{r}^{s}=\begin{cases}B_{r_2+s-1}^{s-1} &\text{if $r=0$}\\
\widetilde{e}^{\epsilon_r(B^{s}_{r-1})}_rB_{r-1}^{s} &\text{if $1\leq r\leq r_2-1$}\\
\widetilde{e}^{\epsilon_{2r_2+s-r}(B^{s}_{r-1})}_{2r_2+s-r}B_{r-1}^{s} &\text{if $r_2\leq r\leq r_2+s$}\end{cases}
$$
\begin{thm}\label{enonar}
Let $s_1\leq s_2$ and $A\otimes B\in B^{r_1,s_1}\otimes B^{r_2,s_2}$. Then we have
$$H(A\otimes B)=-\sum_{\begin{subarray}{c}1\leq p\leq \mathbf r\\\mathbf{\widetilde{r}}\leq q\leq n\end{subarray}}a_{p,q}+\sum^{n-r_2}_{s=0}\big(\epsilon_{r_2}(A^{s}_{r_2+s-1})-\varphi_{r_2}(B^{s}_{r_2+s-1})\big)_{+}.$$
%=-\sum_{\begin{subarray}{c}1\leq p\leq \mathbf r\\\mathbf{\widetilde{r}}\leq q\leq n\end{subarray}}(a_n^{n-r_2})_{p,q}
\proof
The local energy function is by definition constant on the classical components of $B^{r_1,s_1}\otimes B^{r_2,s_2}$. Hence by 
Proposition~\ref{energy1} it is enough to show that the highest weight element in the classical component of $A\otimes B$ is of the form $\widetilde{A}\otimes 0$, where the sum over all entries of $\widetilde{A}$ is 
$$\sum_{\begin{subarray}{c}1\leq p\leq \mathbf r\\\mathbf{\widetilde{r}}\leq q\leq n\end{subarray}}a_{p,q}-\sum^{n-r_2}_{s=0}\big(\epsilon_{r_2}(A^{s}_{r_2+s-1})-\varphi_{r_2}(B^{s}_{r_2+s-1})\big)_{+}.$$
We assume $r_1\leq r_2$, since the proof for $r_1\geq r_2$ proceeds similarly.
We claim the following:
\begin{cl}
The entries in the first $s+1$ rows of $B^s_{r_2+s}$ are zero, i.e. $(b^{s}_{r_2+s})_{p,q}=0$ for all $1\leq p \leq r_2$, $r_2\leq q\leq r_2+s$. Especially we have $B^{n-r_2}_{n}=0$.
\end{cl}
We proof the claim by induction on $s$ and start with $s=0$. By the definition of the elements $B^{s}_r$ (see \eqref{seq}) we get 
$$B^{0}_{r_2}=\widetilde{e}^{\epsilon_{r_2}(B^0_{r_2-1})}_{r_2}\cdots\widetilde{e}^{\epsilon_2(B^0_1)}_2\widetilde{e}^{\epsilon_1(B^0_0)}_1B.$$
It follows $(b^0_1)_{1,r_2}=0$, since $\epsilon_1(B^0_1)=\epsilon_1(\widetilde{e}^{\epsilon_1(B)}_1B)=0$ and  $(b^0_2)_{2,r_2}=0$ since $\epsilon_2(B^0_2)=\epsilon_2(\widetilde{e}^{\epsilon_2(B^0_1)}_2B^0_1)=0$. The action of  $\widetilde{e}_2$ on $B^{0}_{1}$ preserves the first column which implies $(b^0_2)_{1,r_2}=0$. By repeating the above arguments we obtain for $B^{0}_{j}$
$$(b^0_j)_{1,r_2}=(b^0_j)_{2,r_2}=\cdots=(b^0_j)_{j,r_2}=0,$$
which shows the initial step. 
Now suppose that the first $s$ rows are zero in $B^s_0=B^{s-1}_{r_2+s-1}$, i.e. $(b^{s-1}_{r_2+s-1})_{p,q}=0$ for all $1\leq p \leq r_2$, $r_2\leq q\leq r_2+s-1$. Again by \eqref{seq}
$$B^{s}_{r_2+s}=\widetilde{e}^{\epsilon_{r_2}(B^{s}_{r_2+s-1})}_{r_2}\cdots\widetilde{e}^{\epsilon_{r_2+s-1}(B^{s}_{r_2})}_{r_2+s-1}\widetilde{e}^{\epsilon_{r_2+s}(B^{s}_{r_2-1})}_{r_2+s}\widetilde{e}^{\epsilon_{r_2-1}(B^{s}_{r_2-2})}_{r_2-1}\cdots\widetilde{e}^{\epsilon_2(B^{s}_{1})}_2\widetilde{e}^{\epsilon_1(B^{s}_{0})}_1B^{s}_{0}.$$
We obtain as above $(b^s_1)_{1,r_2+s}=0$ since the first $s$ rows are zero and $\epsilon_{1}(\widetilde{e}^{\epsilon_1(B^{s}_{0})}_1B^{s}_{0})=\epsilon_1(B^{s}_{1})=0$. A similar consideration as above shows therefore
$$(b^s_{r_2-1})_{1,r_2+s}=(b^s_{r_2-1})_{2,r_2+s}=\cdots=(b^s_{r_2-1})_{r_2-1,r_2+s}=0$$
and hence $B^s_{r_2-1}$ is of the form 

\begin{figure}[H]
$$
\begin{array}{|c|c|c|c|c|} \hline
0 & 0& \dots & 0& 0\\ \hline
\vdots & \vdots & \vdots & \vdots & \vdots \\ \hline
0 & 0 & \dots &  0& \color{red}{*} \\ \hline 
* & * & \dots &  *& * \\ \hline
\vdots & \vdots & \vdots & \vdots & \vdots \\ \hline
* & * & \dots &  *& * \\ \hline  
\end{array}
$$
\caption{The element $B^s_{r_2-1}$}
\label{figure1}
\end{figure}

Applying the operators $$\widetilde{e}^{\epsilon_{r_2+1}(B^{s}_{r_2+s-2})}_{r_2+1}\cdots\widetilde{e}^{\epsilon_{r_2+s-1}(B^{s}_{r_2})}_{r_2+s-1}\widetilde{e}^{\epsilon_{r_2+s}(B^{s}_{r_2-1})}_{r_2+s}B^s_{r_2-1}$$
means that we move the entry $(b^s_{r_2-1})_{r_2,r_2+s}$ (red star in Figure~\ref{figure1}) along the last column to the upper right corner. Finally we apply the remaining operator $\widetilde{e}^{\epsilon_{r_2}(B^{s}_{r_2+s-1})}_{r_2}$ and obtain the desired property for $B^s_{r_2+s}$.\par
Now it is easy to see that the element $A\otimes B$ and $A^{n-r_2}_{n}\otimes B^{n-r_2}_{n}=A^{n-r_2}_{n}\otimes 0$ are in the same classical component and thus $H(A\otimes B)=H(A^{n-r_2}_{n}\otimes 0)$. In order to prove the theorem we shall show
$$H(A^{n-r_2}_{n}\otimes 0)=-\sum_{\begin{subarray}{c}1\leq p\leq r_1\\r_2\leq q\leq n\end{subarray}}a_{p,q}+\sum^{n-r_2}_{s=0}\big(\epsilon_{r_2}(A^{s}_{r_2+s-1})-\varphi_{r_2}(B^{s}_{r_2+s-1})\big)_{+}.$$

So consider the element $A^{n-r_2}_{n}\otimes 0$. There is a sequence $(i_1,\dots,i_l)$ such that 
$\widetilde{e}_{i_l}\cdots\widetilde{e}_{i_1}(A^{n-r_2}_{n}\otimes 0)$ is a highest weight element. Suppose that there exists an integer $k$, $1\leq k\leq l$, with $r_2=i_k$ and let $k$ be minimal with this property. The condition
$$\epsilon_{r_2}\big((\widetilde{e}_{i_{k-1}}\cdots\widetilde{e}_{i_1}A^{n-r_2}_{n})\big)\leq s_1\leq s_2=\varphi_{r_2}(0)$$
yields $\widetilde{e}_{i_k}\cdots\widetilde{e}_{i_1}(A^{n-r_2}_{n}\otimes 0)=0$, which is a contradiction. Thus $i_k\neq r_2$ for all $k\in\{1,\dots,l\}$. It means that the sum $\sum_{\begin{subarray}{c}1\leq p\leq r_1\\r_2\leq q\leq n\end{subarray}}(a^{n-r_2}_n)_{p,q}$ is stable under the action with $\widetilde{e}_{i_l}\cdots\widetilde{e}_{i_1}$, because the only possibility to change the sum is to apply the operator $\widetilde{e}_{r_2}$. It follows on the one hand that $\widetilde{e}_{i_l}\cdots\widetilde{e}_{i_1}(A^{n-r_2}_{n}\otimes 0)$ is a highest weight element, i.e all entries above the $r_2$-th row are zero by Lemma~\ref{hwe} and on the other hand the sum over the entries below is preserved. By Proposition~\ref{energy1} we have
$$H(A^{n-r_2}_{n}\otimes 0)=-\sum_{\begin{subarray}{c}1\leq p\leq r_1\\r_2\leq q\leq n\end{subarray}}(a^{n-r_2}_n)_{p,q}.$$
Thus it remains to prove 
\begin{equation}\label{rem}\sum_{\begin{subarray}{c}1\leq p\leq r_1\\r_2\leq q\leq n\end{subarray}}(a^{n-r_2}_n)_{p,q}=\sum_{\begin{subarray}{c}1\leq p\leq r_1\\r_2\leq q\leq n\end{subarray}}a_{p,q}-\sum^{n-r_2}_{s=0}\big(\epsilon_{r_2}(A^{s}_{r_2+s-1})-\varphi_{r_2}(B^{s}_{r_2+s-1})\big)_{+}.\end{equation}

Note that 
\begin{flalign*}A^{n-r_2}_{n}=&\widetilde{e}_{r_2}^{\big(\epsilon_{r_2}(A^{n-r_2}_{n-1})-\varphi_{r_2}(B^{n-r_2}_{n-1})\big)_{+}}\widetilde{e}^{\big(\epsilon_{r_2+1}(A^{n-r_2}_{n-2})-\varphi_{r_2+1}(B^{n-r_2}_{n-2})\big)_{+}}_{r_2+1}\cdots\widetilde{e}^{\big(\epsilon_n(A^{n-r_2}_{r_2-1})-\varphi_n(B^{n-r_2}_{r_2-1})\big)_{+}}_n\times &\\&\times \widetilde{e}^{\big(\epsilon_{r_2-1}(A^{n-r_2}_{r_2-2})-\varphi_{r_2-1}(B^{n-r_2}_{r_2-2})\big)_{+}}_{r_2-1}\cdots\widetilde{e}^{\big(\epsilon_1(A^{n-r_2}_{0})-\varphi_1(B^{n-r_2}_{0})\big)_{+}}_{1}A^{n-(r_2+1)}_{n-1},\end{flalign*}
\begin{flalign*}A^{n-(r_2+1)}_{n-1}=&\widetilde{e}_{r_2}^{\big(\epsilon_{r_2}(A^{n-(r_2+1)}_{n-2})-\varphi_{r_2}(B^{n-(r_2+1)}_{n-2})\big)_{+}}\widetilde{e}^{\big(\epsilon_{r_2+1}(A^{n-(r_2+1)}_{n-3})-\varphi_{r_2+1}(B^{n-(r_2+1)}_{n-3})\big)_{+}}_{r_2+1}\times\cdots\times &\\&\times\widetilde{e}^{\big(\epsilon_{n-1}(A^{n-(r_2+1)}_{r_2-1})-\varphi_{n-1}(B^{n-(r_2+1)}_{r_2-1})\big)_{+}}_{n-1} \widetilde{e}^{\big(\epsilon_{r_2-1}(A^{n-r_2}_{r_2-2})-\varphi_{r_2-1}(B^{n-r_2}_{r_2-2})\big)_{+}}_{r_2-1}\times\cdots\times&\\&\times\widetilde{e}^{\big(\epsilon_1(A^{n-r_2}_{0})-\varphi_1(B^{n-r_2}_{0})\big)_{+}}_{1}A^{n-(r_2+2)}_{n-2}.\end{flalign*}
\vdots
\begin{flalign*}\hspace{-5cm}A^{0}_{r_2}=\widetilde{e}^{\big(\epsilon_{r_2}(A^{0}_{r_2-1})-\varphi_{r_2}(B^{0}_{r_2-1})\big)_{+}}_{r_2}\cdots\widetilde{e}^{\big(\epsilon_2(A^{0}_{1})-\varphi_2(B^{0}_{1})\big)_{+}}_2\widetilde{e}^{\big(\epsilon_1(A^{0}_{0})-\varphi_r(B^{0}_{0})\big)_{+}}_1A.\end{flalign*}
Hence the sum over the entries in the last $(n-r_2+1)$ rows of $A^{n-r_2}_{n}$ is given by \eqref{rem}, since the only possibility to decrease the sum $\sum_{\begin{subarray}{c}1\leq p\leq r_1\\r_2\leq q\leq n\end{subarray}}a_{p,q}$ is to apply $\widetilde{e}_{r_2}$.
\endproof
\end{thm}

\begin{ex}
Let $r_1=r_2=n$. The energy of $A\otimes B$ is given by
\begin{align*}&H(A\otimes B)=\sum^{n}_{j=1}a_{j,n}+&\\&
+\Big(a_{n,n}+\Big(a_{n-1,n}+\cdots+\Big(a_{2,n}+\big(a_{1,n}-\varphi_1(B)\big)_{+}-\varphi_2(B)\Big)_{+}-\cdots-\varphi_{n-1}(B)\Big)_{+}-\varphi_n(B)\Big)_{+}.\end{align*}
\end{ex}
%%%%%%%%%%%%%%%%%%%%%%%%%%%%%%%%%%%%%%%%%%%%%%%%%%%%%%%%%%%%%%%%%%%%%%%%%%%%%%%%%%%%%%%%%%%%%%%%%%%%%%%%%%%%%%%%%%%%%%%%%%%%%%%%%%%
%         
%%%%%%%%%%%%%%%%%%%%%%%%%%%%%%%%%%%%%%%%%%%%%%%%%%%%%%%%%%%%%%%%%%%%%%%%%%%%%%%%%%%%%%%%%%%%%%%%%%%%%%%%%%%%%%%%%%%%%%%%%%%%%%%%%%%
\section{Perfect crystals and the ground-state path}\label{section5}
\subsection{Perfect crystals}
We recall the notion of a perfect crystal first introduced in \cite{KKMMNN1992}. Let $c=\sum^n_{i=1}a^{\vee}_i\alpha^{\vee}_i$ be the canonical central element associated to $\Lg$ and $P^+=\{\Lambda\in P\mid \lambda(\alpha^{\vee}_i)\in\Z_{\geq 0}\}$ be the set of dominant integral weights. The level of $\Lambda\in P^{+}$ is defined as $\lev(\Lambda):=\Lambda(c)$. For $\ell\in\Z_{\geq 0}$ let
$$P^+_{\ell}=\{\Lambda\in P^+\mid \lev(\Lambda)=\ell\}.$$

\begin{defn}\label{pc}
A $\U^{'}_q(\Lg)$-crystal $B$ is called a perfect crystal of level $\ell>0$, if the following conditions are satisfied:
\begin{enumerate}
\item $B$ is isomorphic to the crystal graph of a finite-dimensional $\U^{'}_q(\Lg)$-module.
\item $B\otimes B$ is connected.
\item there exists a classical weight $\lambda_0\in P_0$ such that
$$\wt(B)\subseteq \lambda_0+\sum_{i\neq 0}\Z_{\leq 0}\alpha_i,$$
and there is a unique element in $B$ of weight $\lambda_0$.
\item For any $b\in B$, we have $\lev\big(\sum_{i\in I}\epsilon_i(b)\Lambda_i\big)\geq \ell$
\item For all $\Lambda\in P^+_{\ell}$, there exist unique elements $b_{\Lambda},b^{\Lambda}$, such that
$$\sum_{i\in I}\epsilon_i(b_{\Lambda})\Lambda_i=\Lambda=\sum_{i\in I}\varphi_i(b^{\Lambda})\Lambda_i.$$
\end{enumerate}
\end{defn}
\begin{ex}\mbox{}
\begin{enumerate}
\item Let $\Lg=A^{(1)}_2$, the $\U^{'}_q(\Lg)$-crystal
$$
\begin{tikzpicture}

  \node (n1) at (0,3)  {$A$};
  \node (n2) at (2,3)  {$B$};
  \node (n3) at (4,3)  {$C$};
  
 \draw[->] (n1)edge node[above,midway] {1}(n2);
 \draw[->] (n2)edge node[above,midway] {2}(n3);
 \draw[->] (n3) to[out=-160, in=-30]node[below] {0} (n1);
 
\end{tikzpicture}
$$
is a perfect crystal of level $1$
\item Let $\Lg=C^{(1)}_2$, the $\U^{'}_q(\Lg)$-crystal
$$
\begin{tikzpicture}

  \node (n1) at (0,3)  {$A$};
  \node (n2) at (2,3)  {$B$};
  \node (n3) at (4,3)  {$C$};
  \node (n4) at (6,3)  {$D$};
  
 \draw[->] (n1)edge node[above,midway] {1}(n2);
 \draw[->] (n2)edge node[above,midway] {2}(n3);
 \draw[->] (n3)edge node[above,midway] {1}(n4);
 \draw[->] (n4) to[out=-160, in=-30] node[below] {0} (n1);
 
\end{tikzpicture}
$$
is not a perfect crystal of level 1 because we may have $b_{\Lambda_1}=B \text{ or }D$.
\end{enumerate}
\end{ex}
\subsection{Path realization}
Perfect crystals are of particular importance because one can give a path realization of affine highest weight crystals via perfect crystals. Let $B$ be a perfect crystal of level $\ell$ and $\Lambda=\sum_{i\in I}a_i\Lambda_i$ be a dominant integral weight with $\lev(\Lambda)=\ell$. The crystal graph $B(\Lambda)$ associated to the affine Lie algebra $\Lg$ can be realized as follows.\par

Let 
$$\Lambda_0=\Lambda,\ \Lambda_{k+1}=\sum_{i\in I}\epsilon_i(b^{\Lambda_k})\Lambda_i;\ b_k=b^{\Lambda_k}.$$
The sequence 
$$\mathbf p_{\Lambda}=(b_k)^{\infty}_{k=0}=\cdots b_k\otimes b_{k-1}\otimes\cdots\otimes b_0$$
is called the ground-state path of weight $\Lambda$ and a sequence 
$$\mathbf p=(p_k)^{\infty}_{k=0}=\cdots p_k\otimes p_{k-1}\otimes\cdots\otimes p_0$$
with the property $p_k\in B$ for all $k$ and $p_k=b_k$ for all $k>>0$ is called a $\Lambda$-path. We have the following important theorem from \cite{KKMMNN92}.
\begin{thm}\label{mthm}
There exists an isomorphism of crystals
$$\Psi: B(\Lambda)\longrightarrow P(\Lambda), u_{\lambda}\mapsto \mathbf p_{\Lambda}$$
where $P(\Lambda)$ is the set of all $\Lambda$-path in $B$.
\end{thm}
Due to a result of \cite{FOS10} the Kirillov-Reshetikhin crystal $B^{i,\ell}$ is perfect for all non-exceptional types if $\ell$ is a multiple of $c_i$ (see Figure 4 in \cite{ST12}). We give an alternative proof for the perfectness for $A^{(1)}_{n}$ by using the model from \cite{K12} and describe for any $\ell$ and $\Lambda\in P_{\ell}^+$ the elements $b_{\Lambda},b^{\Lambda}\in B^{i,\ell}$.
\begin{thm}\label{perfectness}
The Kirillov-Reshetikhin crystal $B^{i,\ell}$ is perfect of level $\ell$ and for $\Lambda=\sum_{i\in I}a_i\Lambda_i\in P^+_{\ell}$ we have
$$b^{\Lambda}=
\begin{array}{|c|c|c|c|c|} \hline
a_{i+1} & a_{i+2}& \dots & a_{2i}& a_{2i+1}\\[4pt] \hline
a_{i+2} & a_{i+3} & \dots &a_{2i+1}& a_{2i+2} \\[4pt]\hline
\vdots & \vdots & \vdots & \vdots & \vdots \\ \hline
a_{n-i+1} & a_{n-i+2} & \dots &a_{n-1}& a_{n} \\[4pt]\hline
a_{n-i+2} & a_{n-i+3} & \dots &a_{n}& a_{0} \\[4pt]\hline
\vdots & \vdots & \vdots & \vdots & \vdots \\ \hline
a_{n} & a_{0} & \dots &a_{i-3}& a_{i-2} \\[4pt]\hline
a_{0} & a_{1} & \dots &  a_{i-2}& a_{i-1} \\[4pt] \hline 
   \end{array} 
\quad
b_{\Lambda}=
\begin{array}{|c|c|c|c|c|} \hline
a_{1} & a_{2}& \dots & a_{i-i}& a_{i}\\[4pt] \hline
a_{2} & a_{3} & \dots &a_{i}& a_{i+1} \\[4pt]\hline
\vdots & \vdots & \vdots & \vdots & \vdots \\ \hline
a_{n-i} & a_{n-i+1} & \dots &a_{n-2}& a_{n-1} \\[4pt]\hline
a_{n-i+1} & a_{n-i+2} & \dots &  a_{n-1}& a_{n} \\[4pt] \hline 
   \end{array}
$$
\proof
We check stepwise the properties from Definition~\ref{pc}, where (1) and (3) are obviously true.  
In order to prove that $B^{i,\ell}\otimes B^{i,\ell}$ is connected we show that any highest weight element in the tensor product is connected to $0\otimes 0$, i.e. from each highest weight element there is a path in the crystal graph of $B^{i,\ell}\otimes B^{i,\ell}$ to $0\otimes 0$. So let $A\otimes 0$ be an arbitrary highest weight element. We prove the claim by induction over the sum over the entries of $A$. If $A=0$, there is nothing to prove. So let $j$ be maximal
with $a_{i-j,i+j}\neq 0$. Then we consider the element 
$$\widetilde{A}\otimes 0=\widetilde{f}^{a_{i-j,i+j}}_0\widetilde{e}^{a_{i-j,i+j}}_1\cdots\widetilde{e}^{a_{i-j,i+j}}_{i-j-1}\widetilde{e}^{a_{i-j,i+j}}_n\cdots\widetilde{e}^{a_{i-j,i+j}}_{i+j+1}(A\otimes 0),$$
which is a highest weight element and has by induction the desired property. Hence $A\otimes 0$ has the desired property and  condition (2) is proven.
Now let $A\in B^{i,\ell}$ be an arbitrary element. Note
\begin{equation}\label{11}
\epsilon_l(A)=\begin{cases}
  a_{i,i},  & \text{ if $l=i$}\\
  \sum^{i}_{j=q^l_+(A)}a_{j,l}-\sum^{i}_{j=q^l_{+}(A)+1}a_{j,l-1}, & \text{ if $l>i$}\\
  \sum^{q^l_{-}(A)}_{j=i}a_{l,j}-\sum^{q^l_{-}(A)-1}_{j=i}a_{l+1,j}, & \text{ if $l<i$,}
\end{cases}
\geq \begin{cases}
  a_{i,i},  & \text{ if $l=i$}\\
  a_{i,l}, & \text{ if $l>i$}\\
  \sum^{n}_{j=i}a_{l,j}-\sum^{n-1}_{j=i}a_{l+1,j}, & \text{ if $l<i$,}
\end{cases}
\end{equation}
Thus 
\begin{equation}\label{12}\sum^n_{j=0}\epsilon_j(A)=(\ell-\sum^n_{j=i}a_{1,j}-\sum^n_{j=2}a_{j,n})+\sum^n_{j=1}\epsilon_j(A)\geq (\ell-\sum^n_{j=i}a_{i,j})+\sum^n_{j=i}\epsilon_j(A)\geq \ell,\end{equation}
which shows property (4). 
The last property can be deduced from the following observation. Let $b_{\Lambda}$ as above and let $\lev(\Lambda)=\sum_{j\in I} a_j=\ell$, then
$$\epsilon_0(b_{\Lambda})=\ell-\sum_{j\neq 0} a_j=a_0,\ \epsilon_i(b_{\Lambda})=a_i.$$ 
Moreover, for the $j$-th and $(j+1)$-th column ($j<i$) we get $q^j_{-}(b_{\Lambda})=i$ and thus $\epsilon_j(A)=a_j$ and
for the $(j-1)$ and $j$-th row ($j>i$) of $b_{\Lambda}$ we get $q^j_{+}(b_{\Lambda})=i$ and thus $\epsilon_j(A)=a_j$. A similar calculation for $b^{\Lambda}$ proves that both $b^{\Lambda}$ and $b_{\Lambda}$ have the desired property. It remains to prove that $b_{\Lambda}$ is unique (the uniqueness for $b^{\Lambda}$ proceeds similarly). Let $A=(a_{p,q})$ be another element with $\sum_{j\in I}\epsilon_j(A)\Lambda_i=\Lambda$. Since $\sum_{j\in I}\epsilon_j(A)=\ell$ we get with \eqref{11} and \eqref{12}
$$\epsilon_l(A)=\begin{cases}
  a_{i,i},  & \text{ if $l=i$}\\
  a_{i,l}, & \text{ if $l>i$}\\
  \sum^{n}_{j=i}a_{l,j}-\sum^{n-1}_{j=i}a_{l+1,j}, & \text{ if $l<i$,}
\end{cases},
$$
which forces 
$$ p^l_{-}(A)=n\ \forall l=1,\dots,i-1\text{ and }q^l_{+}(A)=i\ \forall l=i+1,\dots,n$$
It follows $\epsilon_i(A)=a_{i,i}=a_i, \epsilon_{i+1}(A)=a_{i,i+1}=a_{i+1},\dots,\epsilon_n(A)=a_{i,n}=a_n$ and therefore the last column of $A$ is the same as the last column of $b_{\Lambda}$. Our aim is to prove that the $(i-1)$-th column is also the same.
Since $q^n_{+}(A)=i$ we get $a_{i-1,n}\leq a_{n-1}$ and since $p^{i-1}_{-}(A)=n$ we get $a_{i-1,n}\geq a_{n-1}$. By repeating this argument with the remaining inidices we obtain $a_{i-1,n}=a_{n-1}, a_{i-1,n-1}=a_{n-2},\dots, a_{i-1,i+1}=a_{i}$. Moreover, 
$\epsilon_{i-1}(A)=\sum^{n}_{j=i}a_{i-1,j}-\sum^{n-1}_{j=i}a_{i,j}=a_{i-1,i}=a_{i-1}$. Consequently the $(i-1)$-th column of $A$ and $b_{\Lambda}$ coincide. By repeating the same method with the $(i-2)$-th and $(i-1)$-th column we obtain that the $(i-2)$-th column of $A$ is the same as the $(i-2)$-th column of $b_{\Lambda}$. We repeat this procedure until we get $A=b_{\Lambda}$.
\endproof
\end{thm}
Finally we can describe easily the ground-state path in $B^{i,\ell}$. We identify $\Lambda=\sum_{j\in I}a_j\Lambda_j$ with the tuple $(a_0,a_1,\dots,a_n)$. Then 
\begin{equation}\label{g}\Lambda_0=(a_0,a_1,\dots,a_n),\ \Lambda_1=(a_i,a_{i+1},\dots,a_n,a_0,\dots,a_{i-1}),\dots\end{equation}
It means that $\Lambda_{k+1}$ arises from $\Lambda_k$ by cutting the tuple $\Lambda_k$ at the $i$-th position into two tuples $(a_0,\dots,a_{i-1})$ and $(a_i,\dots,a_{n})$ and gluing both pieces in reverse order.

$$
\begin{tikzpicture}

  \node (n1) at (-0.8,3)  {$\Lambda_k=(a_0,a_1,\dots,a_n)$};
  \node (n2) at (3,3)  {$(a_0,\dots,a_{i-1})$};
  \node (n3) at (5.5,3)  {$(a_i,\dots,a_{n})$};
  \node (n4) at (10.5,3)  {$\Lambda_{k+1}=(a_i,a_{i+1},\dots,a_n,a_0,\dots,a_{i-1})$};
  
  \draw[->] (n1) to (n2);
   \draw[->] (n3) to(n4); 
  
 \draw[<->,red] (n2) to[out=-30, in=-160] node[below]{}  (n3);
 \draw[<->,red] (n3) to[out=160, in=30] node[above] {}(n2);
  
\end{tikzpicture}$$

The ground state path is then given by 
$$\mathbf p_{\Lambda}=(b_k)^{\infty}_{k=0}=\cdots b_k\otimes b_{k-1}\otimes\cdots\otimes b_0,$$
where $b_k=b^{\Lambda_k}$ with $\Lambda_k$ as in \eqref{g} and $b^{\Lambda_k}$ as in Theorem~\ref{perfectness}.
\begin{ex}\mbox{}
\begin{enumerate}
\item If $i=n$, then 
$$\Lambda_0=(a_0,\dots,a_n),\ \Lambda_1=(a_n,a_0,\dots,a_{n-1}),\ \dots, \Lambda_n=(a_1,\dots,a_n,a_0),\ \Lambda_{n+1}=\Lambda_0$$
and
$$
b_0=
\begin{array}{|c|c|c|c|} \hline
a_{0} & a_{1}& \dots & a_{n-1}\\ \hline
\end{array}
,\
b_1=
\begin{array}{|c|c|c|c|} \hline
a_{n} & a_{0}& \dots & a_{n-2}\\ \hline
\end{array}
,\dots,
b_n=
\begin{array}{|c|c|c|c|} \hline
a_{1} & a_{2}& \dots & a_{n}\\ \hline
\end{array}
$$
\item For $i=n-1$ and $n$ is odd, say $n=2j-1$ we get
$$\Lambda_0=(a_0,\dots,a_n),\ \Lambda_k=(a_{n-(2k-1)},a_{n-(2k-2)},\dots,a_n,a_0,\dots,a_{n-2k}),\ 1\leq k\leq j-1 $$

$$b_0=
\begin{array}{|c|c|c|c|c|} \hline
a_{n} & a_{0} & \dots &a_{n-4}& a_{n-3} \\ \hline
a_{0} & a_{1} & \dots &  a_{n-3}& a_{n-2} \\ \hline 
   \end{array}
	,\
	b_k=
\begin{array}{|c|c|c|c|c|c|c|} \hline
a_{n-2k} & a_{n-(2k-1)} & \dots &a_0 &\dots& a_{n-2k-4}& a_{n-2k-3} \\ \hline
a_{n-(2k-1)} & a_{n-(2k-2)} & \dots &a_1 &\dots&  a_{n-2k-3}& a_{n-2k-2} \\ \hline 
   \end{array}
	$$
\end{enumerate}
\end{ex}
% References
%%%%%%%%%%%%%%%%%%%%%%%%%%%%%%%%%%%%%%%%%%%%%%%%%%%%%%%%%%%%%%%%%%%
\bibliographystyle{plain}
\bibliography{kr-crystals-energy-biblist}
\end{document}